\documentclass[11pt]{article}
\usepackage{xy}
\usepackage{amsmath,amsthm}
\usepackage{amstext, latexsym, amssymb,latexsym}
\usepackage{mathrsfs}
\newcounter{question}[section]
\newtheorem{thm}{Theorem}

\def\R{\mathbb R}
\def\C{\mathbb C}

\def\phi{\varphi}
\def\ol{\overline}

\def\lim{\mathop{\rm lim}}

\def\pf{{\it Proof:}~}
\def\qed{$\,_{\square}\\$}

\newtheorem{rem}[thm]{Remark}

\newtheorem{obs}[thm]{Observation}
\numberwithin{equation}{section}

\begin{document}
\baselineskip 7.2 truemm

\title{\textbf {\Large\center Branch points of area-minimizing
projective planes}}

\normalsize
\author{Robert Gulliver}
\date{}
\maketitle \baselineskip=1.2\normalbaselineskip \noindent

\footnotetext{MSC Classification 53A10, 58E12, 49Q05}

\centerline{\it Dedicated to the memory of Robert Osserman}

\begin{abstract}
Minimal surfaces in a Riemannian manifold
$M^n$ are surfaces which are stationary for area: the first
variation of area vanishes. In this paper we focus on surfaces of
the topological type of the real projective plane $\R P^2$. 
We show that a
minimal surface $f:\R P^2\to M^3$ which has the smallest area,
among those mappings which are not homotopic to a constant
mapping, is an immersion. That is, $f$ is free of branch points.
As a major step toward treating minimal surfaces of the type of the
projective plane, we extend the fundamental theorem of branched
immersions to the nonorientable case. We also resolve a question
on the directions of branch lines posed by Courant in 1950.
\end{abstract}

\maketitle

\section{Introduction}
Let $M$ be an $n$-dimensional  Riemannian manifold, and let
$\Sigma$ be a compact surface with boundary which
carries a conformal structure; we do {\bf not} assume
$\Sigma$ is orientable.  Many existence theorems for minimal
surfaces in the literature \cite{D31}, \cite{R}, \cite{C} find
solutions by minimizing the {\bf  energy} of a mapping 
$f:\Sigma \to M$, where both $f$ and the conformal structure of
$\Sigma$ are allowed to vary. The energy may be written as
\begin{equation} \label{energy}
 E(f) := \frac12 \int_\Sigma (|f_x|^2 + |f_y|^2)\, dx\, dy 
\end{equation}
where $(x,y)$ are local conformal coordinates for $\Sigma$, and
subscripts are used to denote partial derivatives. Write
$\frac{D}{\partial x}$, etc., for
covariant partial derivatives in the Riemannian manifold $M$.
If the mapping $f$ and the conformal structure on $\Sigma$ are
stationary for $E$, then the  resulting mapping is {\bf harmonic}: 
\begin{equation} \label{harm}
\Delta f :=
\frac{D}{\partial x}\frac{\partial f}{\partial x} +
\frac{D}{\partial y}\frac{\partial f}{\partial y} =0
\end{equation}
and {\bf conformal}:
\begin{equation} \label{conf}
 |f_x| \equiv |f_y|, \langle f_x, f_y \rangle \equiv 0. 
\end{equation}

We shall refer to a conformally parameterized harmonic mapping as
a {\em conformally parameterized minimal surface} (CMS).  Observe that
for any $W^{1,2}$ mapping, $E(f)$ is bounded below by the area 
%
\begin{equation} \label{area}
 A(f) := \int_\Sigma |f_x \wedge f_y |\, dx\, dy,
\end{equation}
with equality if and only if $f$ is a conformal mapping almost
everywhere.  In particular, a mapping which minimizes $E(f)$ in a
geometrically defined class of mappings has minimum area among
mappings of the admissible class (see \cite{D39}, p. 232 or Remark
\ref{serrin} below).  Moreover, a conformal mapping which is
harmonic, that is, stationary for $E$, is minimal, that is,
stationary for area.

A CMS is an immersion except at a discrete set of {\bf branch
points.} Let a point of $\Sigma$ be given by $(0,0)$ in some local
conformal coordinates $(x,y)$ for $\Sigma$. Write $z=x+iy$. Then
$(0,0)$ is a {\em branch point} of $f$ of {\em order} $m-1$ if for some
system of coordinates $u^1, \dots, u^n$ for $M$ and for some
complex vector $c$, $f(x,y)$ satisfies the asymptotic description 
$$f^1(x,y) + i f^2(x,y) = c z^m + O(z^{m+1})$$ and 
$$ f^k(x,y) = O(z^{m+1}),$$ 
$k=3,\dots, n$, as $(x,y)\to (0,0)$.  Here we have written 
$f^k(x,y)$ for the value of the $k^{th}$ coordinate $u^k$ at 
$f(x,y)$, $k=1, \dots , n$, and $O(z^{m+1})$ denotes any
``remainder"
function bounded by a constant times $|z^{m+1}|$. We shall
refer to a mapping which is an immersion except at a discrete set
of branch points as a {\bf branched immersion} (see \cite{GOR}).

Our main theorem is 
%
%
\begin{thm} \label{main}
Suppose $\Sigma$ is of the topological type of the real projective
plane $\R P^2$. Let a CMS $f:\Sigma^2 \to M^3$ have minimum area
among all $h:\Sigma \to M$ which are not homotopic to a constant
mapping.  Then $f$ is an immersion.
\end{thm}
In order to prove this theorem, we will distinguish between two
types of branch points: see \cite{O}, \cite{GOR}, \cite{Alt72},
\cite{Alt73}, and \cite{G73}. A {\bf false branch point} of a
branched immersion $f:\Sigma \to M$ is a branch point $z_0$
such that the image set $f(U)$ is an embedded surface, under
another parameterization, for some neighborhood $U$ of $z_0$ in
$\Sigma$.  Otherwise, we call it a {\bf true branch point}. A
branched immersion $f:\Sigma \to M$ is said to be {\bf ramified}
if there are two disjoint open sets $V,W \subset \Sigma$ with
$f(V) = f(W)$. If $f$ is ramified in every neighborhood of a point
$z_0$, we say that $z_0$ is a {\bf ramified branch point}.  Note
that any false branch point of a branched immersion must be
ramified.  Osserman showed that in codimension one, a branched
immersion $f:\Sigma^2 \to M^3$ with a true branch point cannot
minimize area, see \cite{O} and Theorem \ref{oss} below,  in
contradiction to assertions of Douglas (p. 239 of \cite{D32}) and
of Courant (footnote p. 46 of \cite{C41}). 
On the other hand, regarding false branch points, we
shall extend to nonorientable surfaces the fundamental theorem of
branched immersions in \cite{G75}, and show that if a branched CMS
$f:\R P^2 \to M^n$ is ramified, with any codimension, then there
is another CMS $\widetilde{f}:\R P^2 \to M$ with at most half the
area of $f$. 

We would like to acknowledge the interest of Simon Brendle in this
problem, whose questions, not used in \cite{BBEN}, stimulated us
to investigate this research topic. We are also indebted to the
late Jim Serrin for pointing us toward Remark \ref{serrin}.

%
%
\section{Analysis of branch points}\label{anal}

This section reports on material that has appeared in the
literature, see especially \cite{G73}. In this paper, we shall
discuss certain steps in the interest of clarity and completeness. 

Let $\Sigma^2$ be a compact surface with a conformal structure, $M^n$ a
Riemannian manifold, and let $f: \Sigma \to M$ be a CMS. Consider a
branch point $z_0 \in \Sigma$ for $f$. Write $D$ for the
Riemannian connection on $M$. Let local conformal
coordinates $(x,y)$ for $\Sigma$ and local coordinates 
$(q_1, \dots, q_n)$ for $M$ be introduced with $z_0 = (0,0)$ and
$f(z_0) = (0,\dots,0)$ in these coordinates. Then equation
\eqref{harm} may be rewritten
$$\frac{D}{\partial\overline{z}}\frac{\partial f}{\partial z}=0,$$
where we write the complex coordinate $z=x+iy$, 
$\frac{\partial f}{\partial z}= 
\frac12[\frac{\partial f}{\partial x}-i\frac{\partial f}{\partial
y}]$, and 
$\frac{D}{\partial \overline z}= 
\frac12[\frac{D}{\partial x} +i \frac{D}{\partial y}].$
In this form,
we see that harmonicity implies that the complex tangent vector
$\frac{\partial f}{\partial z}$ is holomorphic to first order. It
is readily shown that for some positive integer $m$ and for some
complex tangent vector $c=a+ib$ to $M$ at $f(z_0)$, 
%
\begin{equation}
f(z)=\mathcal{R}\{c z^m\} +O_2(|z|^{m+1}).
\end{equation}
Here we write $\mathcal{R}\{v\}$ for the real part of a complex
vector $v$, and we have used the big-O notation with the subscript
$2$, meaning that as $z \to 0$, the remainder term is bounded by a
constant times $|z|^{m+1}$, its first partial derivatives are
bounded by a constant times $|z|^m$ and its second partial
derivatives are bounded by a constant times $|z|^{m-1}$. It
follows from the conformality condition \eqref{conf} that the
complex-bilinear inner product 
$\langle c,c \rangle = |a|^2 - |b|^2 +2i \langle a,b \rangle =0.$ 
Choose a new system of coordinates $p_1, \dots, p_n$ for $M$ near 
$f(z_0)$ with $\frac{\partial}{\partial p_1}=a$ and 
$\frac{\partial}{\partial p_2}=b$; and a new system of coordinates
$(\widetilde{x},\widetilde{y})$ for $\Sigma$ with
$\widetilde{z}=\widetilde{x}+i\widetilde{y}=|a|^{1/m} z$.
Then along the mapping $f$,
$$p_1+ip_2=\widetilde{z}^m +\sigma(\widetilde{z})$$
 and
$$p_\ell=\psi_\ell(\widetilde{z}),$$
$\ell = 3,\dots,n$, where 
$\sigma(\widetilde{z}), \psi_\ell(\widetilde{z}) = O_2(\widetilde{z}^{m+1})$. 
We now define a non-conformal complex parameter $w=u_1+iu_2$ on a
neighborhood of the branch point in $\Sigma:$
%
\begin{equation}\label{defw}
w:=\widetilde{z}\,
\Big[1+\widetilde{z}^{-m}\sigma(\widetilde{z})\Big]^{1/m}.
\end{equation}
Then $w$ is a $C^{1,\alpha}$ coordinate on $\Sigma,$ for some
H\"{o}lder exponent $\alpha >0$, in terms of
which the coordinate representation of $f$ is simplified:
%
\begin{equation}\label{anonpar}
p_1+ip_2=w^m, 
\end{equation} 
and
$$p_\ell=\phi_\ell(w)=O_2(w^{m+1}),$$
$\ell=3,\dots,n.$

We now turn our attention to the case n=3 of {\bf codimension
one}. The self-intersection of the surface is determined by the
single real-valued function $\phi(w)=\phi_3(w).$ Define 
$\overline\phi(w)= \phi(\zeta_m w)$, where 
$\zeta_m = e^{2\pi i/m}$ is a primitive $m^{\it th}$ root of unity,
and let $\Phi(w)=\phi(w)-\overline\phi(w).$ Then the zeroes of $\Phi$
correspond to curves of intersection of the surface with itself.
But both $\phi$ and $\overline\phi$ satisfy the {\it same}
quasilinear minimal surface equation in $M$, with the {\it same}
coefficients. Therefore, their difference 
$\Phi:= \phi-\overline\phi$ satisfies a linear homogeneous PDE:
%
\begin{equation}\label{Phieq}
\sum_{i,j=1}^2 a_{ij}\Phi_{u_i u_j}+\sum_{i=1}^2 a_i \Phi_{u_i} +
a \Phi =0, 
\end{equation}
whose coefficients, as functions of $w$, are obtained by 
integrating from the PDE
satisfied by $\phi$ to the PDE satisfied by $\overline\phi$ along
convex combinations. We have $a_{ij}(0,0) =\delta_{ij}$. 

It follows that $\Phi$ satisfies an asymptotic formula 
%
\begin{equation}\label{asymp}
\Phi(w) = {\mathcal R}\{A w^N\} + O_2(w^{N+1}),
\end{equation}
for some integer $N >m$ and some complex constant $A\neq 0$ 
(see \cite{HW}). We shall call
$N-1$ the {\bf proper index} of the branch point.
We may sketch the proof of Hartman and Wintner in \cite{HW}. We 
rewrite the PDE \eqref{Phieq} in terms of the complex 
gradient $\Phi_w= \frac12(\Phi_{u_1}- i \Phi_{u_2})$. 
If $\nabla\Phi(w)=o(w^{k-1})$, then  we test against the function
$g(w) = w^{-k} (w-\zeta)^{-1}$, where $\zeta \neq 0$ is small, to
show that $\Phi_\zeta(\zeta) = a \zeta^k+o(\zeta^{k+1})$ for some
$a\in \C$. Proceeding by induction on $k$, one finds formula 
\eqref{asymp} for some integer $N$ and some $A\neq 0$, unless 
$\nabla\Phi(w) \equiv 0$. Details are as in \cite{HW}, pp. 455-458.

Thus, there are two alternatives: $\Phi$ is either identically
zero or satisfies the asymptotic formula \eqref{asymp} where 
the integer $N$ is $\geq m+1$ and the complex number
$A\neq 0$. If $\Phi \equiv 0$, then $z_0$ is a false branch point:
see section \ref{false} below.

If $\Phi$ is not $\equiv 0$, then we have a {\bf true branch point}. 

%
%
\begin{thm} \label{oss}
(\cite{O}, \cite{G73}.)  Suppose $\Sigma$ is a 
surface with a conformal structure, $M^3$ a Riemannian manifold 
and $f:\Sigma \to M$ a
mapping which has smallest area in a $C^0$ neighborhood of
$f$. Then $f$ has no true branch points.
\end{thm}
\pf
As we have just seen, a true branch point $z_0$ has an order 
$m-1 \geq 1$ and a coordinate neighborhood in $\Sigma$ with a
$C^{1,\alpha}$ complex coordinate $w$, $w=0$ at $z_0$, such that
$f$ has the representation \eqref{anonpar} near $z_0$. Adjacent
sheets $p_3=\phi(w)$ and $p_3=\overline\phi(w)=\phi(\zeta_m w)$ 
intersect when $\Phi(w)=0$, which, according to formula
\eqref{asymp}, occurs along $2N\geq 6$ arcs in $\Sigma$ forming
equal angles $\pi/N$ when they leave $w=0$. Let one of these
arcs be parameterized as $\gamma_1:[0,\varepsilon]\to \Sigma$, and
let the corresponding arc be $\gamma_2:[0,\varepsilon]\to \Sigma$, 
defined by $w(\gamma_2(t)) =\zeta_m w(\gamma_1(t))$.
Then for all $0\leq t \leq \varepsilon$, 
$\phi(\gamma_1(t))=\phi(\gamma_2(t))$. Note from formula
\eqref{anonpar} that all three coordinates coincide: the mapping
$f(\gamma_1(t))= f(\gamma_2(t)), 0\leq t \leq \varepsilon$.

We may now construct a Lipschitz-continuous and piecewise smooth
surface $\widetilde{f}$ which has the same area as $f$, but has
discontinuous tangent planes, following Osserman \cite{O}. The
idea of the following construction is that the parameter domain
$D$ may be cut along the arcs $\gamma_1((0,\varepsilon))$ and
$\gamma_2((0,\varepsilon))$, opened up to form a lozenge, with two
pairs of adjacent sides originally identified, and then closed up
along the remaing two pairs of adjacent sides. 

In detail: choose an open topological disk $D\subset \C$ on which the
coordinate $w$ is defined, $(0,0)\in D$, and which is invariant 
under the rotation taking $w$ to $\zeta_m w$.  Assume
$\gamma_1(\varepsilon)$ and $\gamma_2(\varepsilon)$ are the first
points along $\gamma_1$ resp. $\gamma_2$ which lie on the boundary of
$D$. We shall construct a discontinuous,
piecewise $C^1$ mapping $Q:B_1 \to D$, such that
$\widetilde{f}(\zeta):=f(Q(\zeta))$ is nonetheless continuous, and
$Q$ is one-to-one and onto except for sets of measure $0$. Here,
$B_1$ is the disk $\{z\in \C: |z|<1\}$. Choose
points $A_i=\gamma_i(\varepsilon/2)$, $i=1,2$. Then $D$ is
broken along $\gamma_1$ and $\gamma_2$ into two curvilinear pentagons 
with vertices $\gamma_1(\varepsilon), A_1, (0,0), A_2$ and
$\gamma_2(\varepsilon).$ The edges of these pentagons are 
$\gamma_1([\varepsilon/2,\varepsilon])$,  $\gamma_1([0,\varepsilon/2])$,
$\gamma_2([0,\varepsilon/2])$, $\gamma_2([\varepsilon/2,\varepsilon])$ 
and one of the two arcs of $\partial D$ with endpoints 
$\gamma_1(\varepsilon)$ and $\gamma_2(\varepsilon)$.
Similarly, break the unit disc $B_1$ 
along the interval $(-1,1)$ of the $x$-axis and the interval
$[-\frac12, \frac12]$ of the $y$-axis into two pentagons.
Each pentagon in $B_1$ will be bounded by four line segments, an
interval along the $y$-axis being used twice, plus the upper or
lower half-circle of $\partial B_1$.  Denote the points
$a=(0,\frac12),$ $e=(0,-\frac12)$, $c_1=(1,0)$, $c_2=(-1,0)$ and
give the origin $(0,0)$ four different names: $b_1$ when approached
from the first quadrant $\{x>0,y>0\}$, 
$b_2$ when approached from the second quadrant $\{x<0, y>0\}$,
$d_2$ when approached from the third quadrant $\{x<0, y<0\}$, and
$d_1$ when approached from the fourth quadrant $\{x>0, y<0\}$. $Q$
will map the pentagon in $B_1$ in the upper half-plane $y>0$ to the pentagon
in $D$ lying counterclockwise from $\gamma_1$ and clockwise from
$\gamma_2$, with $Q(c_1)=\gamma_1(\varepsilon)$, $Q(b_1)=A_1$,
$Q(a)=(0,0)$, $Q(b_2)=A_2$, and $Q(c_2)=\gamma_2(\varepsilon)$.
This describes $Q$ on the boundary of one of the two pentagons;
the other pentagon is similar. The
interior of each pentagon may be made to correspond by a $C^1$
diffeomorphism. We require $Q$ to be continuous along the
$x$-axis. Of course, $Q$ is discontinuous along the intervals
$0<y<\frac12$ and $-\frac12<y<0$
of the $y$-axis. However, we can regain the continuity of
$\widetilde{f}$ by requiring that for $0\leq y \leq \frac12$, along
the interval from $b_2$ to $a$, $Q(0,y)=\gamma_2(\varepsilon y)$ and
along the interval from $b_1$ to $a$, 
$Q(0,y)=\gamma_1(\varepsilon y)$. Similarly, $Q$ will map the
pentagon in the lower half-plane to the pentagon lying
counterclockwise from $\gamma_2$ and clockwise from 
$\gamma_1$, taking care that for $-\frac12 \leq y \leq 0$, 
$Q(0,y)=\gamma_2(-\varepsilon y)$ as approached from the third
quadrant, and $Q(0,y)=\gamma_1(-\varepsilon y)$ as approached from
the fourth quadrant.

Then the continuous, piecewise $C^1$ image surfaces
$\widetilde{f}(B_1)$ and $f(D)$ consist of the same pieces of
surface in $M$, and therefore have the same area. But the tangent
planes of $\widetilde{f}$ are discontinuous along the $y$-axis at
each $y$ in the open intervals $-\frac12 < y < 0$ and 
$0 < y < \frac12$, which implies that $\widetilde{f}$, and
therefore $f$, does not have minimum area: smoothing
$\widetilde{f}$ near these arcs reduces its area to first order.

Note that by choosing $\varepsilon$ small, we may make $\widetilde{f}$
arbitrarily close to $f$ in the uniform topology.
\qed

\begin{obs}\label{courant} 
R. Courant (see p. 123 of \cite{C}) asked whether the curves of 
self-intersection near any true branch point of a minimal surface in
Euclidean $\R^3$ meet at equal angles at the branch point. We may
observe that this conjecture is partially correct, as seen above:
the curves of intersection given by $\Phi(w)=0$ make equal angles
at the branch point.  However, the totality of the curves of
intersection may 
make a variety of angles at a branch point of order 
$m-1\geq 3$. Specifically, if $m=4$ and the proper index $N-1=5$,
then the zeroes of $\Phi(w)=\phi(w)-\phi(\zeta_4 w)$ will occur
along $2N=12$ curves forming equal angles in the parameter plane,
and therefore (since angles in the $w$-plane at $0$ are multiplied
by $m=4$ in $\R^3$) equal angles in the
tangent plane to $f(\Sigma)$ at the branch point.  But their images are
only those curves of intersection coming from {\bf successive} pairs of
the four ``sheets" of the surface. We also have the intersection
of {\bf non-successive} sheets, which are the images of zeroes of
$\Phi_2(w):=\phi(w)-\phi(-w)$, since $\zeta_4^2=i^2=-1$. The leading
term $A w^6$ of the asymptotic formula \eqref{asymp} cancels when
we compute $\Phi_2(w)=\phi(w)-\phi(-w)=\Phi(w)+\Phi(iw)$. 

It follows from the PDE satisfied by $\phi(w)$ and by $\phi(-w)$ 
that $\Phi_2$ satisfies an
elliptic PDE analogous to \eqref{Phieq}, and therefore an
asymptotic relation analogous to \eqref{asymp} with leading term
$A_2w^{N_2}$ for some integer $N_2>N$ and some 
$A_2\in \C\backslash \{0\}$. That is, the curves of intersection
of non-successive sheets form a family of equally spaced
directions, which are presumably independent of the directions of
the curves of intersection of successive sheets.  This philosophy
is justified by the following explicit example with $N=6$ and
$N_2=7$.

Choose $a,b\in\C\backslash \{0\}$. Using the Weierstra\ss\ representation
(see \cite{O69}, p. 63) for a minimal surface $f:\C \to \R^3$ in
Euclidean 3-space, based on the polynomials $4z^3$ and
$2az^2+2bz^3$ (the latter representing the Gau\ss\ map in
stereographic projection), we have the specific CMS with 
\begin{eqnarray}
f^1_z(z) = \Big[1-(az^2+bz^3)^2\Big]2z^3\\
f^2_z(z) = -i\Big[1+(az^2+bz^3)^2\Big]2z^3\\
f^3_z(z) = 4z^3(az^2+bz^3),
\end{eqnarray}
which leads to 
$$w^4:=f^1+if^2 = z^4-\frac{\ol{a}^2\ol{z}^8}{2}-
\frac{8}{9}\ol{a}\ol{b}\ol{z}^9- \frac{2}{5}\ol{b}^2\ol{z}^{10}$$
and to
$$z=w\Big(1+\frac{\ol{a}^2\ol{w}^8}{8 w^4}+
\frac{2\ol{a}\ol{b}\ol{w}^9}{9 w^4}+O_2(|w|^6)\Big),$$
via an extensive, but straightforward, computation. 
Recall that each component $f^k$ of
$f$ is real and harmonic as a function of $z$. Rewriting 
$$f^3(z)=8\,\mathcal{R} \{ \frac{a}{6} z^6 + \frac{b}{7} z^7\}$$
as a (non-harmonic) function of $w$, we find 
$$ \phi(w) = 8\,\mathcal{R}\{\frac{a}{6}w^6 + \frac{b}{7}w^7+
\frac{\ol{a}|a|^2}{8} \ol{w}^8 w^2 +O_2(|w|^{11})\}:$$
the difference of $\phi$ on successive sheets is 
\begin{equation}\label{succ}
\Phi(w):=\phi(w)-\phi(iw)=
8\,\mathcal{R}\{\frac{a}{3}w^6+\frac{b}{7}(1+i)w^7\}+ O_2(|w|^{10})
\end{equation}
and on non-successive sheets is
\begin{equation}\label{nonsucc}
\Phi_2(w):=\phi(w)-\phi(-w)=
8\,\mathcal{R}\{bw^7\}+ O_2(|w|^{10}).
\end{equation}

From the formula \eqref{succ}, we see that $N=6$, $A=\frac{8a}{3}$ 
and the curves of intersection of
successive sheets are curves in $\R^3$ leaving the branch point
along the $(x_1,x_2)$-plane, which is the tangent plane to
$\Sigma$ at the branch point, in the $12$ directions
$(\cos(4\theta),\sin(4\theta),0)$, where $6\theta+\arg(a)$ is an
integer multiple of $\pi$. The $12$ directions are
paired off to form $6$ curves in $\R^3$ leaving the branch point 
at equal angles $\frac{2\pi}{3}$.

Similarly, from the formula \eqref{nonsucc}, we see that $N_2=7$,
$A_2=8b$ and the curves 
of intersection of nonsuccessive sheets are seven curves leaving
the branch point and making equal angles (images of 14 curves in
the $w$-plane, paired). The arguments of 
$a,b \in \C\backslash\{0\}$ may be given arbitrary values, so that the
angle between a representative of the family of six curves of
self-intersection and a representative of the family of seven
curves of self-intersection may be chosen arbitrarily. For most
choices, these 13 curves in $\R^3$ will {\bf not} form equal
angles at the branch point.

\end{obs}

%
%
\section{False branch points}\label{false}

The elimination of false branch points from an area-minimizing CMS
$f:\Sigma \to M$ is in general only possible by comparison with
surfaces $\Sigma_0$ of {\it reduced topological type} (see
\cite{D39}, \cite{G77}): for orientable surfaces, $\Sigma_0$ has 
smaller genus or the same total genus and more connected components. 
As an oriented example, we may
choose $\Sigma$ to be a surface of genus $2$ and $\Sigma_0$ to be
a torus. Then there is a branched covering $\pi: \Sigma \to
\Sigma_0$ with two branch points of order one. (Think of $\Sigma$
as embedded in $\R^3$ so that it is invariant under a rotation by
$\pi$ about the $z$-axis, and meets the $z$-axis only at two
points: the quotient under this rotation is a torus.) Now choose a
minimizing CMS $f_0: \Sigma_0 \to M^n$, and let $f:\Sigma \to M^n$
be $f = f_0 \circ \pi$. Then $f$ has two false branch points.  In
order to be sure that $f$ minimizes area in its homotopy class, we
may choose $M^3$ to be a flat $3$-torus with two small periods and
one large period.  As one sees from this example, in order to show 
that false branch points do not
occur, we must assume that $f$ minimizes area among mappings from
surfaces of the topological  type of $\Sigma$ {\bf and} of lower
topological type. This hypothesis was used by J. Douglas (see
\cite{D39}), in a strict form, to find the existence of minimal
surfaces $f:\Sigma \to \R^3$ with prescribed boundary. For 
$\R P^2$, however, there are no nonorientable surfaces of lower type.

Results in the literature for false branch points have until now
assumed that $\Sigma$ is {\it oriented}, see \cite{G73},
\cite{Alt73}, \cite{GOR}, \cite{G75}, \cite{G77} and \cite{T}.
In order to treat false branch points for nonorientable
surfaces, we will need to extend certain known results. In
particular, the following theorem appears in \cite{G75} for 
{\em orientable} surfaces, possibly with boundary, including
surfaces of prescribed mean curvature vector not necessarily zero.

%
%
\begin{thm}\label{fund}
(Fundamental theorem of branched immersions)
Let $\Sigma^2$ be a compact surface with boundary endowed with a 
conformal structure, $\partial\Sigma$ possibly empty 
and $\Sigma$ not necessarily orientable.  Let $M^n$ be a Riemannian
manifold and $f:\Sigma \to M$ a CMS. Assume that the restriction
of $f$ to $\partial\Sigma$ is injective. Then there exists a compact
Riemann surface with boundary $\widetilde\Sigma$, a branched covering 
$\pi:\Sigma \to \widetilde\Sigma$ and a CMS 
$\widetilde{f}: \widetilde\Sigma \to M$ such that 
$f = \widetilde{f} \circ \pi$.  Moreover, the restriction of
$\widetilde{f}$ to $\partial\widetilde\Sigma$ is injective.
Further, $\widetilde\Sigma$ is orientable if and only if $\Sigma$
is orientable.
\end{thm}
\pf
If $\Sigma$ is orientable, then 
Theorem 4.5 of \cite{G75} provides an orientable quotient surface
$\widetilde\Sigma$, a branched covering 
$\pi:\Sigma \to \widetilde\Sigma$ and an unramified CMS 
$\widetilde f: \widetilde\Sigma \to M$ such that 
$f=\widetilde f \circ \pi$. 

There remains the case where $\Sigma$ is {\bf not orientable}.
Assume, without loss of generality, that $\Sigma$ is connected. 

Let $p:\widehat\Sigma \to \Sigma$ be the oriented double cover of
$\Sigma$, with the induced conformal structure. Then 
$\widehat\Sigma$ is connected and orientable, and $p$ is two-to-one. 
The composition 
$\widehat{f}=f \circ p: \widehat\Sigma \to M$ is a CMS, defined on
an orientable surface, and we may apply Theorem 4.5 of \cite{G75}
to find a compact orientable surface with boundary 
$\widehat{\widetilde\Sigma}$, an unramified
CMS $\widehat{\widetilde{f}}:\widehat{\widetilde\Sigma}\to M$ and
an orientation-preserving branched covering 
$\widehat{\pi}:\widehat\Sigma \to \widehat{\widetilde\Sigma}$ so that
$\widehat{f}$ factors as $\widehat{\widetilde{f}}\circ\widehat{\pi}$.

Now let $\widetilde\Sigma$ be the quotient surface of
$\widehat{\widetilde\Sigma}$ under the identification of 
$\widehat{\pi}(x^+) \in \widehat{\widetilde\Sigma}$ with 
$\widehat{\pi}(x^-) \in \widehat{\widetilde\Sigma}$ whenever 
$x^\pm \in \widehat{\Sigma}$ and $p(x^+)= p(x^-)$ in $\Sigma$.
Then for each $x \in \Sigma$, $p^{-1}(x)$ consists of two points
$x^+,\ x^-\in\widehat{\Sigma}$ and there are diffeomeorphic
neighborhoods of $\widehat{\pi}(x^+)$ and of $\widehat{\pi}(x^-)$
which are thereby identified in $\widetilde\Sigma$, with reversal
of orientation. This implies 
that $\widetilde\Sigma$ is a differentiable $2$-manifold.  Write
$\widetilde{p}:\widehat{\widetilde\Sigma}\to\widetilde\Sigma$ for
the quotient mapping. Then $\widetilde{f}:\widetilde\Sigma \to M$
is well defined such that 
$\widehat{\widetilde{f}} = \widetilde{f}\circ\widetilde{p}$.
Also, for $x\in \Sigma$, the two pre-images $x^+, x^-\in
\widehat{\Sigma}$ have
$\widetilde{p}\circ\widehat{\pi}(x^+)=\widetilde{p}\circ\widehat{\pi}(x^-),$ 
so that we may define
$\pi:\Sigma \to \widetilde\Sigma$ by 
$\pi(x):= \widetilde{p}\circ \widehat{\pi}(x^\pm)$.

Note that the mappings $p, \widehat{\pi}$ and $\widetilde p$ are
surjective, and therefore also $\pi:\Sigma\to\widetilde{\Sigma}$.

In the event that $\partial{\Sigma}$ is nonempty, since
$f=\widetilde f \circ \pi$ restricted to $\partial\Sigma$  is
injective, it follows readily that the restriction of
$\widetilde{f}$ to $\partial{\widetilde{\Sigma}}$ is injective. 

Then in the above construction, for each $x\in \Sigma$, 
$f$ defines the
same piece of surface, with opposite orientations, on
neighborhoods of $x^+$ and of $x^-$. The branched covering
$\widehat\pi: \widehat\Sigma \to \widehat{\widetilde\Sigma}$
preserves orientation, implying that
$\widehat\pi(x^+) \neq \widehat\pi(x^-)$.
Since $\widehat{\widetilde\Sigma}$ is connected, there is a path
from $\widehat\pi(x^+)$  to $\widehat\pi(x^-)$ whose image in
$\widetilde\Sigma$ reverses orientation. Therefore
$\widetilde\Sigma$ is not orientable. 
\qed

%
%
\section{An immersion of $\R P^2$}

We are now ready to give the proof of the main Theorem \ref{main}.
Let $f: \R P^2 \to M^3$ be a CMS into a three-dimensional
Riemannian manifold, which has minimum area among all mappings
$\R P^2 \to M^3$ not homotopic to a constant. Write
$\Sigma = \R P^2$. From Theorem \ref{oss}, we see that $f$ has no
true branch points. (For this conclusion, it would suffice that $f$
minimizes area in a $C^0$ neighborhood of each branch point.)

There remains the possibility of false branch points.

We first recall the computation of the Euler characteristic of a
surface. For a compact, connected surface which is either
orientable or nonorientable, the Euler characteristic
$$ \chi(\Sigma) = 2 - r(\Sigma),$$
where $r(\Sigma)$ is the topological characteristic of $\Sigma$
\cite{D39}, also known as the nonorientable genus; 
we
shall adopt the term {\bf demigenus}. If $\Sigma$ is orientable,
then it has even demigenus and genus $\frac12 r(\Sigma)$. If 
it is non-orientable, then $\Sigma$ may be constructed
by adding $r(\Sigma)$ cross-caps to the sphere.
The demigenus of the sphere equals zero, of $\R P^2$ equals one, of
the torus and the Klein bottle equals two. For other compact
surfaces without boundary, the demigenus is $\geq 3$. 

Now according to Theorem \ref{fund}, there is a compact Riemann
surface $\widetilde\Sigma$, a branched covering
$\pi:\Sigma \to \widetilde\Sigma$ and an unramified CMS
$\widetilde f:\widetilde\Sigma \to M$ such that
$f = \widetilde f \circ \pi$. We will apply the Riemann-Hurwitz
formula to the branched covering $\pi$:
%
\begin{equation}\label{RH}
 \chi(\Sigma) = d \, \chi(\widetilde\Sigma) - {\mathcal O}(\pi),
\end{equation}
where $d$ is the degree of $\pi$, 
${\mathcal O}(\pi)$ is the total order of branching of $\pi$, and
$\chi$ is the Euler number.  Suppose that $f$ has a false
branch point, or more generally, a ramified branch point.  Then
${\mathcal O}(\pi)\geq 1,$ and the branched covering $\pi$ has
degree $d \geq 2$.

Using the formula \eqref{RH}, we can determine the
topological type of $\widetilde\Sigma$. Since $\Sigma$ is
homeomorphic to $\R P^2$, it has $\chi(\Sigma)=1$.  We also know
that $d>0$
and ${\mathcal O}(\pi) \geq 1$. It follows that the demigenus
$r(\widetilde\Sigma) \leq 1$. Otherwise, the integer
$r(\widetilde\Sigma)$ would be $\geq 2$, which implies
$\chi(\widetilde\Sigma) \leq 0$ and by the formula \eqref{RH}, 
$1 \leq -{\mathcal O}(\pi)\leq -1$,  a contradiction.  That is,
$\widetilde\Sigma$ is either the sphere or the projective plane.
But according to Theorem \ref{fund}, since $\Sigma$ is not
orientable, $\widetilde\Sigma$ is not orientable; therefore,
$\widetilde\Sigma$ is homeomorphic to $\R P^2$.

Note that if $\widetilde f:\widetilde\Sigma \to M$ were homotopic
to a constant mapping, then so would be $f = \widetilde f \circ
\pi$.

On the other hand,
the area of $f$ equals the area of $\widetilde f$ times the degree
$d$ of $\pi$. But $d\geq 2$, so the area of $\widetilde f$ is at
most one-half the area of $f$.  But this would mean that $f$ does
not have minimum area among maps $:\R P^2 \to M$ not homotopic to
a constant mapping, contradicting our hypothesis.  This implies
that $f$ has no branch points, and is therefore an immersion.  
\qed

\begin{rem}
We have treated conformally parameterized {\em minimal} surfaces
in this paper. However, the proofs go through with only minor
changes for projective planes of nonzero {\em prescribed} mean
curvature $H:M^3 \to \R$, provided that $f(\Sigma)$ has a
transverse orientation.  This can occur only when $M$ is
non-orientable.

It also appears plausible that a version of Theorem \ref{fund} can be
extended to the more general case of mappings satisfying the {\em
unique continuation property}, see \cite{GOR}.  
\end{rem}

\begin{rem}\label{serrin}
An alternative approach to branch points of minimal surfaces
of the type of the disk appears in the recent book \cite{T} by 
Tromba. The second and higher variations of energy $E$ 
(see \eqref{energy}) of a CMS $f$ are computed in a
neighborhood of a branch point $z_0$, and the lowest nonvanishing
variation is shown to be negative if the branch point is
nonexceptional. This is defined in terms of the {\bf index} $i>m$ 
of $z_0$, where $i+1$ is the order of contact of the mapping
with the tangent plane at the branch point. Note that the proper
index $N$ is $\geq i$ (recall the definition \eqref{asymp}).
If $i+1$ is an integer multiple of $m$, where $m-1$ is the order
of the branch point, then the branch point is called {\em
exceptional}. Tromba also shows that exceptional interior branch
points will not occur, provided that the mapping has minimum area
$A$ among surfaces with the same boundary curve. 

In fact, minimizing area and minimizing energy, under such Plateau
boundary conditions, are equivalent properties of a Lipschitz
continuous mapping $f:B \to \R^n$, where $B$ is the unit disk in
$\R^2$.  Namely, if $E(f) \leq E(g)$ for all Lipschitz-continuous
mappings $g:B \to \R^n$ defining the same boundary curve, then
$A(f) \leq A(g)$ for all such $g$, as we now show. 

Otherwise, for some $g:B \to \R^n$ with the same Plateau boundary 
conditions as $f$, $A(f) > A(g)$. Write $\eta = A(f)-A(g)>0$.
Approximate $g$ with $g^\delta(w):=(g(w),\delta w) \in \R^{n+2}$.
Then $g^\delta$ is a Lipschitz immersion, so there are conformal
coordinates $\widetilde w=:F^{-1}(w)$ for some bi-Lipschitz 
homeomorphism $F:B \to B$ which preserves $\partial B$ (see
\cite{M}). Write
$\widetilde g^\delta(\widetilde w):= g^\delta(F(\widetilde w))$ for
the conformal mapping with the same image as $g^\delta$. Define
$\widetilde g:B \to \R^n$ by composing $\widetilde g^\delta$ with
the projection from $\R^{n+2} \to \R^n$.  Then the energy 
$E(\widetilde g^\delta) = E(\widetilde g) + \delta^2 E(F)$.
Also, the area $A(g^\delta) \leq A(g) + C\delta$ for some constant
$C$. It follows that
$$E(\widetilde g)\leq E(\widetilde g^\delta)=A(\widetilde
g^\delta)=A(g^\delta)\leq A(g) + C\delta <A(f) \leq E(f)$$
if $\delta$ is chosen small enough that $C\delta <\eta$. This
implies that $f$ does not minimize energy, a contradiction.

The converse implication, that if $f$ minimizes area then it
minimizes energy, follows similarly.

We may observe that with this remark, Tromba's book \cite{T} gives
an independent proof that a mapping from the disk into $\R^3$
which minimizes energy for prescribed Plateau boundary conditions
is an immersion in the interior. 

In addition, the book contains new, partial results on
boundary branch points.

\end{rem}

\vspace{1cm}
\ \  \\
\noindent Robert Gulliver\\
School of Mathematics, University of Minnesota,\\
Minneapolis, MN 55455, USA\\
{\tt e-mail:gulliver@math.umn.edu}\\


\begin{thebibliography}{123}
\bibitem{Alt72} H. Wilhelm Alt, Verzweigungspunkte von H-Fl\"{a}chen
I, Math. Z. {\bf 127} (1972), 333--362.
\bibitem{Alt73} H. Wilhelm Alt, Verzweigungspunkte von H-Fl\"{a}chen
II, Math. Annalen {\bf 201} (1973), 33--56.
\bibitem{BBEN}  Hugh Bray, Simon Brendle, Michael Eichmair and Andre Neves, 
Area-minimizing projective planes in 3-manifolds. Comm. Pure Appl.
Math. {\bf 63} (2010), 1237–1247.
%
\bibitem{C41} Richard Courant, On a generalized form of Plateau's
problem, Trans. Amer. Math. Soca {\bf 50}, 40--47(1941).
%
\bibitem{C} Richard Courant, Dirichlet's principle, Conformal
Mapping and Minimal Surfaces. New York: Wiley Interscience 1950.
\bibitem{D31} Jesse Douglas, Solution of the problem of Plateau,
Trans. Amer. Math. Soc. {\bf 33} (1931), 263--321. 
\bibitem{D32} Jesse Douglas, One-sided minimal surfaces with a
given boundary, Trans. Amer. Math. Soc. {\bf 34} (1932), 731-756.
\bibitem{D39} Jesse Douglas, Minimal Surfaces of Higher
Topological Structure, Annals of Math. {\bf 40} (1939), 205--298.
%
\bibitem{G73} Robert Gulliver, Regularity of minimizing surfaces
of prescribed mean curvature, Annals of Math. {\bf 97} (1973),
275--305.
\bibitem{G75} Robert Gulliver, Branched immersions of surfaces and
reduction of topological type, I, Math. Z. {\bf 145} (1975),
267--288.
\bibitem{G77}Robert Gulliver, Branched immersions of surfaces and
reduction of topological type, II, Math. Ann. {\bf 230} (1977),
25--48.
\bibitem{GOR} Robert Gulliver, Robert Osserman and Halsey Royden,
A theory of branched immersions of surfaces, Amer. J. Math {\bf
95} (1973), 750--812.
%
\bibitem{HW} Philip Hartman and Aurel Wintner, On the local
behavior of solutions of non-parabolic partial differential
equations, Amer. J. Math. {\bf 75} (1953), 258--287.
\bibitem{M} Charles B. Morrey, On the solutions of quasi-linear
elliptic partial differential equations.
Trans. Amer. Math. Soc. {\bf 43} (1938), 126–166. 
\bibitem{O} Robert Osserman, A proof of the regularity everywhere
of the classical solution to Plateau's problem, Annals of Math.
{\bf 91} (1970), 550--569.
\bibitem{O69} Robert Osserman, A Survey of Minimal Surfaces, New
York: Van Nostrand Reinhold 1969.
\bibitem{R} Tibor Rad\'{o}, The problem of least area and the
problem of Plateau, Math. Z. {\bf 32} (1930), 763--796.
\bibitem{T} Anthony Tromba, A theory of branched minimal surfaces.
Springer Verlag 2012.

\end{thebibliography}
\end{document}